\documentstyle{amsart}
\title[Mckay-Thompson Series]{The mckay-Thompson
 series associated with the
 irreducible characters of the monster}

\author{Koichiro Harada}
\author{Mong Lung Lang}

\begin{document}

\baselineskip=12pt


\maketitle

\begin{abstract}  Let $ \Bbb V = \coprod_{h = 0}^{\infty} \Bbb V_h$
 be the graded monster module of the monster simple group $\Bbb M$
 and let $\chi_k$
 be an irreducible representation of $\Bbb M$. The generating function
 of  $c_{hk}$ (the multiplicity of $\chi_k$ in $\Bbb V_h$) is determined.
  Furthermore, the invariance group of  the modular function
  associated with the generating function is also determined in this paper.
\end{abstract}

\section{Introduction}

Let $\Bbb M$ be the monster simple group and $\Bbb V$ be the monster module
 of Frenkel-Lepowsky-Meurman [3]. $\Bbb V$ is a graded $\Bbb M$ module
 $$ \Bbb V = \coprod_{h = 0}^{\infty} \Bbb V_h$$
such that
$$j(q) - 744 = \sum_{ h=0}^{\infty} \text{dim} \Bbb V_h q^{h-1}.$$
In particular, $\text{dim}\Bbb V_0 = 1$, $\text{dim}\Bbb V_1 = 0$, $\text{dim}
\Bbb V_2 = 196884$,
 $\text{dim}\Bbb V_3 = 21493760$, $ \cdot \cdot \cdot \, $.
   Let $\chi_k,$ $1 \le k \le 194$, be the
irreducible characters of $\Bbb M$, which will  often be used to
 denote the irreducible representations also. For the first few
 $\Bbb V_i$'s, we have the decompositions :
$$\Bbb V_0 = \chi_1,$$
$$\Bbb  V_2 = \chi_1 + \chi_2,$$
$$\Bbb V_3 = \chi_1 + \chi_2 + \chi_3.$$
In general, write
$$\Bbb V_h = \sum_{k = 1}^{194}c_{hk}\chi_k$$
 where $c_{hk}$ is the multiplicity of $\chi_k$ in $\Bbb V_h.$
  The table of $c_{hk}$ for $0 \le h \le 51, $ $ 1 \le k \le 194$
 can be found in McKay-Strauss [6].

We also list some of the multiplicities $c_{h1}$ of the trivial
 character $\chi_1$ in $V_h$.
{\small
$$\begin{array}{rrrrrrrrrrrrrrrr}
h & 0 & 2 &3 & 4 &5 & 6 & 7 & 8 & 9 & 10 &11
\, \cdot \cdot \cdot  &
20 \, \cdot \cdot \cdot &
30 \,\cdot \cdot \cdot & 40\,  \cdots & \, 50  \\
c_{h1}&1 &1 & 1 & 2& 2 & 4 &4 &  7& 8 & 12 &
14 \, \cdot \cdot \cdot &
167\,  \cdot \cdot \cdot &
1762\, \cdot \cdot \cdot &15913 \, \cdots & \, 129734\\
\end{array}  $$}
Let us consider, for each irreducible character $\chi_k$, the
 generating function :
$$t_{\chi_k}(x) = x^{-1} \sum_{h = 0}^{\infty} c_{hk} x^h.$$
  The multiplicity $c_{hk}$ can be computed as follows :
$$c_{hk} = \frac{1}{|\Bbb M|} \sum_{g \in \Bbb M} Tr(g|\Bbb V_h)\chi_k(g).$$
Therefore the generating function of $c_{hk}$ is
$$t_{\chi_k}(x) =x^{-1} \sum_{h= 1}^{\infty} \sum_{g \in \Bbb M}
 \frac{1}{|\Bbb M|} Tr(g|\Bbb V_h)\chi_k(g) x^h.$$
If we replace the indeterminate $x$ by $q = e^{2 \pi i z}$, $z \in \Bbb H$
 =  $\{ z \in \Bbb C | \text{Im}(z) > 0 \}$ then
$$t_{\chi_k}(q) = \frac{1}{|\Bbb M|} \sum_{g \in \Bbb M} \chi_k(g) t_g(q)$$
where
$$t_g(q) = q^{-1} \sum_{h=0}^{\infty} Tr(g|\Bbb V_h)q^h$$
is the McKay-Thompson series for the element $g$
 in $\Bbb M$.  Thus $t_{\chi}(q)$ for the irreducible character $\chi$ is
the weighted sum of the McKay-Thompson series for the element $g$ of $\Bbb M$.
  Not all $t_{\chi}(q)$'s are distinct and in fact there are exactly
 172 distinct $t_{\chi}(q)$'s, since
 $$t_{\chi}(q) = t_{\bar \chi}(q)$$
where $\bar \chi$ is the complex conjugate of $\chi$ and there is no
 other equalities among $t_{\chi}(q)$'s.
  One of the obvious questions one will raise here will be :

\vspace{.2cm}

{\bf Problem.}  Determine the invariance group $\Gamma_{\chi}$
  of $t_{\chi}(q).$

\vspace{.2cm}
Here $\Gamma_{\chi}$ is defined to be :

\vspace{.2cm}

{\bf Definition.}  $\Gamma_{\chi} = \{ A \in  SL_2(\Bbb R) | t_{\chi}(Az)
 = t_{\chi}(z)\}.$

\vspace{.2cm}

Since $t_{\chi}(z)$ is a modular function, $\Gamma_{\chi}$ is
 a discrete subgroup of $SL_2(\Bbb R)$.  Let us here review some of the
 properties of the  invariance group $\Gamma_{g}$ of the
 McKay-Thompson series $t_g(z)$ for the element $g \in \Bbb M$.
\vspace{.2cm}

(0).  For $G \subset GL_2^{+}(\Bbb R)$, $\bar G$ is the image of $G$
  in $PGL_2^{+}(\Bbb R)$.

(1).  $\Gamma_0(N) = \{ A =\left( \begin{array}{ll}
                               a & b \\
                               c & d  \\
                          \end{array} \right) \in SL_2(\Bbb Z)
                            | c \equiv 0$ (mod $N$) $\}.$

(2).  For an exact divisor $e||N$ (i.e. $e|N$ and \text{gcd}$(e, \frac{N}{e})
=1$) let
$$W_e = \left( \begin{array}{ll}
                ae &  b \\
               cN  &  de \\
        \end{array}
     \right ),\quad  a, \, b, \, c, \, d \in \Bbb Z, \quad ade^2 - bcN = e. $$
Then $\bar W_e$ normalizes $\bar \Gamma_0(N)$ and
 $\bar W_e^{2} \in \bar \Gamma_0(N).$

(3).  Let $h$ be a divisor of $n$.  Then
 $n|h+ e,f, \cdot \cdot \cdot $  is defined to be
 $$\left(\begin{array}{ll}
                              \frac{1}{h} & 0 \\
                                   0 & 1 \\
                              \end{array}
                           \right)
 \langle \Gamma_0(\frac{n}{h}), W_e,W_f, \cdot \cdot \cdot \rangle
   \left(\begin{array}{ll}
                                  h & 0 \\
                                   0 & 1 \\
                              \end{array}
                           \right).$$

(4).  For each $g$ in $\Bbb M$, $\Gamma_g$,
 the invariance group of $t_g(z),$
 is a normal subgroup
 of index $h_g$ in $n_g|h_g + e_g, f_g, \cdot \cdot \cdot $,
 the eigen group of $g$ [2]. Note that for each $A$ in
  $n_g|h_g + e_g, f_g, \cdot \cdot \cdot$
, $(t_g|A)(z)$ = $\sigma t_g(z)$
 where $\sigma$ is an $h_g$-th root of unity.
We will often use $n$, $h$, $e$, $f$, etc. instead of $n_g$, $h_g$
 $e_g,$ $f_g$, etc.
 for simplicity.

(5).  $\Gamma_g$ contains $\Gamma_0(n_gh_g).$

\vspace{.2cm}

For each irreducible character of $\Bbb M$, we now define :

\vspace{.2cm}

{\bf Definition.}  $N_{\chi}$ = lcm$\{n_gh_g | g \in\Bbb M,
  \chi(g) \ne 0\}.$

\vspace{.2cm}

It is obvious that
  $t_{\chi}(z)$ is invariant under $\Gamma_0(N_{\chi})$.
 Note that  $N_{\chi}$
 can be quite large
($N_{\chi_1}$ = $2^6 3^3 5^2 7\cdot11 \cdot 13\cdot
 17\cdot 19\cdot 23\cdot 29\cdot 31\cdot 41\cdot 47\cdot 59 \cdot 71$)
  or relatively small ($N_{\chi_{166}}$ = $2^63^27 = 4032$).

\vspace{.2cm}

The purpose of this paper is to show

\vspace{.2cm}

{\bf Theorem.}  $\Gamma_{\chi} = \Gamma_0(N_{\chi}).$

\section{Poles of $t_{\chi}(z)$}

For each cusp $c$ in $\Bbb Q \cup \{\infty\}$,
  we define  $\Phi_c$  to be the set
  $\Phi_c$= $\{ g \in \Bbb M | c$ is equivalent to $\infty$ in $\Gamma_g \}$
   and decompose $t_{\chi}(z)$ into :
 $$ \frac{1}{|\Bbb M|} \sum_{g \in \Phi_c}
  \chi(g)t_g(z)
  + \frac{1}{|\Bbb M|} \sum_{g \notin \Phi_c}
  \chi(g)t_g(z). $$
 Since the McKay-Thompson series $t_g(z)$ is a generator
 of the function field  of the compact Riemann surface $
   \Gamma_g \setminus \Bbb H^*$ ($\Bbb H^* = \Bbb H \cup \{\infty \}$)
 of genus 0 and has a unique pole at $\infty$ (and at all
 cusps  $c \in \Bbb Q$ equivalent to $\infty$ in $\Gamma_g$).
    Obviously, $ \frac{1}{|\Bbb M|} \sum_{g \notin \Phi_c}
  \chi(g)t_g(z) $
  is holomorphic
at $c$.
 Hence,
  whether   $c$ is a pole of $t_{\chi}(z)$ or not  is
  determined  by  the singular part of
 $$ \frac{1}{|\Bbb M|} \sum_{g \in \Phi_c}
  \chi(g)t_g(z) $$
 at $c$.
   For example, if $c = \infty$, then
   $\infty$ is a pole of $t_{\chi}(z)$
 if and only if $\chi$ is the trivial character since the singular
 part of  $t_{\chi}(z)$ at $\infty$ is given by
 $\frac{1}{|\Bbb M|} \sum_{ g \in \Bbb M}\chi(g)\frac{1}{q}$ and
$$ \sum_{ g \in \Bbb M}\chi_i(g) = \cases
   |\Bbb M| & \text{if $i = 1$} \\
    0  & \text{if $i \not= 1$} \\
    \endcases $$

Suppose $c \in \Bbb Q$. For each $g$ in $\Phi_c$, let
 $ A \in SL_2(\Bbb Z)$ be chosen
 so that $A\infty = c$. Then $t_g(Az) = (t_g|A)(z)$ has an expansion
 in $q= e^{2 \pi i z}$ of the form
$$ aq^{-\frac{1}{\mu}} + \cdot \cdot \cdot$$
where $\mu = [ \langle \left(\begin{array}{ll}
                       1  &  1\\
                       0 &1   \\
                    \end{array}
                   \right)  \rangle : A^{-1} \Gamma_0(n_gh_g)_cA],$
 where the subscript $c$ denotes the stabilizer.
  We contend that the contribution of the $t_g(z)$ to the singular part
 of $t_{\chi}(z)$ is
$$aq^{-\frac{1}{\mu}}.$$
Indeed, by our assumption the cusp $c$ is equivalent to $\infty$ in
$\Gamma_g$ and so there is $ B \in \Gamma_g$ such that $B\infty = c$
 and
$$(t_g|B)(q) = t_g(q) =  q^{-1} + \sum_{i \ge 0} a_iq^i.$$
 The only difference between $(t_g|A)(z)$ and $(t_g|B)(z)$ lies in the
power of $q$ and $a$, hence the contention.

In order to determine whether $c$ is pole of $t_{\chi}(z)$ or not,
  one has to :

(1). Determine whether
  $c$ is equivalent to $\infty$  in $\Gamma_g$ or not.

(2). Determine the singular part of
   $ \frac{1}{|\Bbb M|} \sum_{g \in \Phi_c}
  \chi(g)t_g(z)$ at $c$.

We will investigate those questions
 in the next section.

\bigskip

\section{equivalence of cusps}

In this section, we will study the equivalence of cusps
 in $\Gamma_g$, $g \in \Bbb M.$

\bigskip

{\bf Lemma 1.} {\em For each exact divisor $e$ of $N$ and
 for  each $c$ such that $\text{gcd}(c,e) = 1$, $\Gamma_0 (N)$  admits
 an Atkin-Lehner involution of the form $W_e = \left( \begin{array}{ll}
                                           ae & b \\
                                           cN & de \\
                                           \end{array} \right).$
 Moreover, one can choose
 either $a$ = 1 or $d$ = 1 if desired. }

\bigskip

{\em Proof}.  For each $c$ such that $c$ and $e$ are relatively prime,
   we have $ \text{gcd}(\frac{cN}{e},e) = 1$.
Hence, there exists $b$ and $y$ such that $ye - \frac{bcN}{e} = 1,$
 or $ye^2 - bcN = e.$ The lemma
follows by writing $ y$ into $ ad $ for suitable
 $a$ and $d$.
\qed

\bigskip

{\bf Lemma 2.} {\em Let $\text{gcd}(a,b) = 1$
 and $M$ be nonzero integers.
  Then there exists a pair of integers $x$ and $y$
 satisfying  $\text{gcd}(xM,y) = 1$
  and $ax + by  = 1.$}

\bigskip

{\em Proof.}  This is a well known fact of the elementary number theory.
 Let $x'$ and $y'$ be a pair of integral solutions
 of  the equation $ax+by=1$ and let $M = M_aM_{y'}M'$ be the
 decomposition of $ M$
 into a product of coprime factors such that
  $M_a$ and $a$,  $M_{y'}$ and $y'$,
  have the same prime factors.  It is clear that
$y = y'+aM'$ and $ x= x'-bM'$ is also a pair of solutions to
 the equation. Note that $\text{gcd}(x,y)=1$ since
 it is a solution of $ax+by=1.$   Furthermore, one has
 $\text{gcd}(y,M)$ = $\text{gcd}(y'+aM',M_aM_{y'}M')$ =
 $\text{gcd}(y'+aM',M')= 1.$
 Therefore  $x$ and $y$ is pair of integral solutions
 of the equation such that $\text{gcd}(xM,y)=1.$
\qed

\bigskip

{\bf Lemma 3.} {\em   Let $\frac{x}{y}$, \text{gcd}($x$,$y$) = 1,
 be a rational number.  Then
$\frac{x}{y}$ is equivalent to some $ \frac{x'}{y'}$, \text{gcd}($x'$,$y'$) = 1
 in $\Gamma_0 (N)$,
where $y'$ = \text{gcd}($N$,$y$).  Furthermore, if $\frac{x}{y}$  is equivalent
 to $\frac{ x''}{y''}$ with $y''|N,$ $y'' > 0$, and
 gcd($x'',y''$) = $1$, then $y'' = y'$.}

\bigskip

{\em Proof.}  Consider the equality
$$\left(
\begin{array}{ll}
a & b \\
cN & d \\
\end{array}
\right) \frac{x}{y} = \frac{ax+by}{cNx+dy}
 = \frac{ax+dy}{y'(cx\frac{N}{y'} + d\frac{y}{y'})} $$
Note that $\text{gcd}(x\frac{N}{y'}, \frac{y}{y'}) = 1 $, hence the
 equation
$$ c\frac{xN}{y'} + d\frac{y}{y'} = 1 $$
 is solvable for $c$, $d$ in  $\Bbb Z$.  Applying Lemma 2,
 we may assume that $c$ and $d$ are integral solutions
of the above  equation such that
$ \text{gcd} ( cN, d) = 1$.
Let $a$ and $b$ be chosen such that
$ ad -cbN = 1$.  Summerizing, we now conclude
 that $\frac{x}{y}$ is equivalent to $\frac{ax+by}{y'}$
by $\pmatrix a&b\\cN&d\endpmatrix \in \Gamma_0(N) $.
  Since $\text{gcd}(ax+by,y')
  = \text{gcd}(ax,y') = \text{gcd}(a,y') = 1$
, first part of the lemma is
 proved.  As for the second part,
suppose
  $$\left(
\begin{array}{ll}
a & b \\
cN & d \\
\end{array}
\right) \frac{x}{y} = \frac{ax+by}{cNx+dy}
 = \frac{ax+dy}{y'(cx\frac{N}{y'} + d\frac{y}{y'})} = \frac{x''}{y''}. $$
 We first note $y'|y''$, since gcd($ax+by,y'$) = $1$. To show $y''|y'$,
 suppose that  $y''$ possesses a prime power $q^t$ such that $y'$ is not
 a multiple of $q^t,$ then  $q^t|y'(cx\frac{N}{y'} + d\frac{y}{y'})$ implies
   $q|(cx\frac{N}{y'} + d\frac{y}{y'})$.  Since $y''|N$,
 $q$ is a divisor of $\frac{N}{y'}$, hence $q|d$. This implies
 that
  $\text{gcd}(cN,d) \ne 1,$
  against our choice of  $c$, $d$.  Thus,  $y''|y'$
  and the second part of the lemma is proved.
\qed

\bigskip

In Lemma 4 and Lemma 5, $ G = n|h + e, f, \cdot \cdot \cdot$  is the
eigen group of the invariance group $\Gamma_g$.

\bigskip

{\bf Lemma 4.} {\em  Let $g \in \Bbb M $ and let $\Gamma_g \le G$ =
 $ n|h + e, f, \cdot \cdot \cdot$
 be the invariance group of $t_g(z)$.
  Then $G\infty$ = $\Gamma_g \infty$}.
\bigskip

{\em Proof.}  Since $ q$ = exp(2$\pi i z$) is a local parameter
of $t_g(z)$, the stabilizer $(\Gamma_{g})_{\infty}$ of $\infty$
is generated by $\left( \begin{array}{ll}
1 & 1 \\
0 &1 \\
\end{array}
\right)$.
As for $G$, $G_{\infty}$
 is generated by
 $\left( \begin{array}{ll}
1 & \frac{1}{h} \\
0 &1 \\
\end{array}
\right)$.
  Hence $$[G : \Gamma_g] = h =
 [ G_{\infty}: (\Gamma_{g})_{\infty}] $$
Consequently, $G \infty = \Gamma_g \infty $.
 \qed

\bigskip

{\bf Lemma 5.} {\em  Let $g \in \Bbb
M $ and let $\Gamma_g \le n|h + e, f, \cdot \cdot \cdot $
be the invariance group of $t_g(z)$. Then $\frac{x}{y}$,
 $\text{gcd}(x,y) = 1$,  is equivalent
to $\infty$ in $\Gamma_g$ if and only if
$$\text{gcd} ( \frac{n}{h}, \frac{y}{\text{gcd}(y,h)} ) \in
   \left\{ \frac{n}{h}, \frac{n}{he}, \frac{n}{hf},
   \cdot \cdot \cdot \right\}$$ }

\bigskip

{\em Proof}.  To simplify our notation, let N = $\frac{n}{h}$.
Choose the Atkin-Lehner involution $W_e$ as  described in
 Lemma 1  with $\text{gcd}(c,N) = 1.$
 One has $W_e \infty$ = $\frac{a}{c\frac{N}{e}},$ and
 $$ \text{gcd}(a, c\frac{N}{e}) = 1.$$
By Lemma 3, there exists $\gamma_e \in \Gamma_0(N)$ such that
 $\gamma_e W_e \infty = \frac{e'}{\frac{N}{e}}$ where
 $ \text{gcd}(e',\frac{N}{e}) = 1 $, since $\text{gcd}(c\frac{N}{e},N) =
 \frac{N}{e}.$
  Note that $e$ is an exact
 divisor of $N$, hence among the representitives
 of  inequivalent cusps
 of $\Gamma_0(N),$  there is exactly one and only one
  cusp $z$ with denominator
 $\frac{N}{e}$ (see Harada [4]).  Without loss of generality, we may assume
that $z$ = $\frac{1}{\frac{N}{e}}$.  Therefore,  we may
 assume that $\gamma_e$ is chosen so that
$\gamma_e W_e \infty = \frac{1}{\frac{N}{e}}$. Hence the
$G$-orbit of $\infty$ can be decomposed into,
$$\left( \begin{array}{ll} \frac{1}{h} & 0 \\
                              0 & 1    \\
                            \end{array} \right) \Gamma_0(N)
 \frac{1}{N}
 \cup
\left( \begin{array}{ll} \frac{1}{h} & 0 \\
                              0 & 1    \\
                            \end{array} \right) \Gamma_0(N)
 \frac{1}{\frac{N}{e}}
 \cup
\left( \begin{array}{ll} \frac{1}{h} & 0 \\
                              0 & 1    \\
                            \end{array} \right) \Gamma_0(N)
 \frac{1}{\frac{N}{f}}
  \cup \cdot \cdot \cdot .$$

Hence $\frac{x}{y}$ is equivalent to $\infty$ in $G$ if
and only if
$$ \left( \begin{array}{ll}
                 h & 0 \\
                 0 &1 \\
          \end{array}
          \right) \frac{x}{y} =
\frac{\frac{hx}{\text{gcd}(y,h)}}{\frac{y}{\text{gcd}(y,h)}} \in
 \Gamma_0(N) \frac{1}{N}
 \cup
 \Gamma_0(N) \frac{1}{\frac{N}{e}}
 \cup
 \Gamma_0(N) \frac{1}{\frac{N}{f}}
  \cup \cdot \cdot \cdot,  $$
which is equivalent to, by Lemma 3,
$$\text{gcd}(\frac{y}{\text{gcd}(y,h)},N) \in  \left\{N,\frac{N}{e},
\frac{N}{f}, \cdot \cdot \cdot
 \right\}.$$
 $G\infty = \Gamma_g\infty$ as shown in Lemma 4 and so
  $\frac{x}{y}$ is equivalent
to $\infty$ in $\Gamma_g$ if and only if
$$\text{gcd}( \frac{y}{\text{gcd}(y,h)}, \frac{n}{h} ) \in
  \left\{ \frac{n}{h}, \frac{n}{he}, \frac{n}{hf},
   \cdot \cdot \cdot \right \}.$$
\qed

\bigskip

{\bf Corollary 6}. {\em  0 = $\frac{0}{1}$ is equivalent to $\infty$
in $\Gamma_g$ if and only if $n = h$ or
$G = n|h + e,f, \cdot \cdot \cdot $ admits
 the Atkin-Lehner involution $W_{\frac{n}{h}}$}.

\bigskip

{\em Proof}.  Since $\text{gcd}(1,\frac{n}{h}) = 1$, $G$ must admits
an Atkin-Lehner involution $W_e$ such that $\frac{n}{he} = 1,$
 hence  $e$ = $\frac{n}{h}$.
\qed

\bigskip

 Let $\chi$ be an irreducible character of the monster $\Bbb M.$
 In order to determine the singular part of
  $$ \frac{1}{|\Bbb M|} \sum_{g \in \Phi_c}
  \chi(g)t_g(z)$$ at the cusp $c = \frac{x}{y}$,
  where  $\text{gcd}(x,y) = 1$ and $y|N_{\chi}$,
   it is necessary to find
 a matrix $P_c$ in $SL_2(\Bbb R)$  such that $P_c\infty = c$
  and determine the $q$-expansion of
 $$ \frac{1}{|\Bbb M|} \sum_{g \in \Phi_c}
  \chi(g)(t_g|P_c)(z),$$
which will be called the $q$-expansion of $t_{\chi}(z)$ at $c.$
Such a matrix $P_c$ is easy to find and  choice is
 not unique.  To ease the computation  of the $q$-expansion of
   $$ \frac{1}{|\Bbb M|} \sum_{g \in \Phi_c}
  \chi(g)(t_g|P_c)(z),$$
 it is necessary to find a  good
  $P_c$ so  that the transformation formula
  of $t_g|P_c$ can be obtained for
 every $g$ in $\Phi_c$ simultaneously.
 What we shall do is as follows. Namely,
 we will  find a matrix  $P_c$ so that one can associate
 with $P_c$ an  upper triangular  matrix ${U_{c,g}}$ such that
  $$P_c{U_{c,g}}^{-1} \in n_g|h_g+e_g, f_g,  \cdot \cdot \cdot \, .$$
Since $n_g|h_g + e_g, f_g, \cdot \cdot \cdot  $ is the
 eigen group of $\Gamma_g$,
 elements in $n_g|h_g + e_g,f_g,  \cdot \cdot \cdot $ map $t_g$ to
 $\sigma_g t_g$ where $\sigma_g$ is an $h_g$-th root of unity
 which depends on $g$ and on some other quantities.  Therefore
 $$t_g|P_c = \sigma_g  t_g|{U_{c,g}}.$$
A good transformation formula for $t_g|P_c$ is obtained
 since $U_{c,g}$ is upper trianglar.

  Let $y_0$ be the exact
 divisor of $N_{\chi}$ such that $y$ and $y_0$ share the
 same prime divisors.
  Then gcd$(y,x\frac{N_{\chi}}{y_0}) = 1$
 and so there is a matrix   $P_c \in SL_2(\Bbb Z)$ of the form
$$ P_c = \left(\begin{array}{ll}
            x &  w  \\
            y & \frac{zN_{\chi}}{y_0}  \\
        \end{array}
        \right).$$
  Lemma 4 implies  that
 $\frac{x}{y}$ is equivalent to $\infty$ in $\Gamma_g$
 if and only if $\frac{x}{y}$ is equivalent
 to $\infty$ in the eigen group $n_g|h_g + e_g, f_g,  \cdot \cdot \cdot$
  of $g$ and so by Lemma 5, $\frac{x}{y}$ is equivalent to $\infty$ in
$\Gamma_g$ if and only if
 $\text{gcd}(\frac{n_g}{h_g}, \frac{y}{\text{gcd}(y,h_g)}) \in
 \left\{ \frac{n}{h}, \frac{n}{he}, \frac{n}{hf},
   \cdot \cdot \cdot \right \}.$
  More precisely, $\frac{x}{y}$ is
 equivalent to $\infty$  by an element in
 $ n_g|h_g$ if
 $$
\text{gcd}(\frac{n_g}{h_g}, \frac{y}{\text{gcd}(y,h_g)}) = \frac{n_g}{h_g}$$
 and   is equivalent to $\infty$ by an
 Atkin-Lehner involution  $W_{e_g}$ of $n_g|h_g + e_g,
 f_g, \cdot \cdot \cdot$ if
 $$   \frac{n_g}{h_ge_g} =\text{gcd}(\frac{n_g}{h_g}, \frac{y}
{\text{gcd}(y,h_g)})
      \in  \left\{ \frac{n}{he}, \frac{n}{hf},
   \cdot \cdot \cdot \right \}.$$

\bigskip

{\bf  Lemma 7.}  {\em  Suppose that
 $\text{gcd}(\frac{n_g}{h_g}, \frac{y}{\text{gcd}(y,h_g)}) = \frac{n_g}{h_g}.$
  Let  $u_g$ be  chosen so that  $ \frac{yu_g}{h_g} +
\frac{zN_{\chi}\text{gcd}(h_g,y)}{y_0h_g}$
    is an integer.
 Then
$$P_cU_{c,g}^{-1} = P_c\left(\begin{array}{ll}
     \frac{h_g}{\text{gcd}(h_g,y)}  &  \frac{u_g}{h_g}   \\
           0  &             \frac{\text{gcd}(h_g,y)}{h_g}  \\
           \end{array}
    \right) \in  n_g|h_g$$
where
$$U_{c,g} = \left(\begin{array}{ll}
     \frac{\text{gcd}(h_g,y)}{h_g}  &  -\frac{u_g}{h_g}   \\

       0  &             \frac{h_g}{\text{gcd}(h_g,y)}  \\
    \end{array}
\right).$$}

\bigskip

{\em Proof.}  To show the existence of $u_g$, simply solve the equation
$$\frac{y}{\text{gcd}(h_g,y)} u_g + \frac{N_{\chi}}{y_0} z \equiv 0
 \pmod {\frac{h_g}{\text{gcd}(h_g,y)}}.$$
 Then $yu_g + z \frac{N_{\chi}}{y_0}  \text{gcd}(h_g,y) \equiv 0 \pmod { h_g}$
 and hence
  $ \frac{yu_g}{h_g} +
\frac{zN_{\chi}\text{gcd}(h_g,y)}{y_0h_g}$
    is an integer.

The matrix
 $$P_c\left(\begin{array}{ll}
     \frac{h_g}{\text{gcd}(h_g,y)}  &  \frac{u_g}{h_g}   \\
           0  &             \frac{\text{gcd}(h_g,y)}{h_g}  \\
           \end{array}
    \right)  =
  \left(\begin{array}{ll}
     \frac{xh_g}{\text{gcd}(h_g,y)}  &   \frac{xu_g+w\text{gcd}(h_g,y)}{h_g}
\\
     \frac{yh_g}{\text{gcd}(h_g,y)} & \frac{yu_g}{h_g} +
                              \frac{zN_{\chi}\text{gcd}(h_g,y)}{y_0h_g} \\
           \end{array} \right)$$
 has the property

(1).   $ \frac{yu_g}{h_g} +
\frac{zN_{\chi}\text{gcd}(h_g,y)}{y_0h_g}$
    is an integer by our choice of $u_g$, and,

(2). $\frac{yh_g}{\text{gcd}(h_g,y)} $ is a multiple of $n_g$ since
$\text{gcd}(\frac{n_g}{h_g}, \frac{y}{\text{gcd}(y,h_g)}) = \frac{n_g}{h_g}$.

Therefore

$$ P_cU_{c,g}^{-1} = P_c\left(\begin{array}{ll}
     \frac{h_g}{\text{gcd}(h_g,y)}  &  \frac{u_g}{h_g}   \\
           0  &             \frac{\text{gcd}(h_g,y)}{h_g}  \\
           \end{array}
    \right) \in  n_g|h_g. \qed $$

\bigskip

{\bf Corollary 8.} {\em  Suppose that
$\text{gcd}(\frac{n_g}{h_g}, \frac{y}{\text{gcd}(y,h_g)}) = \frac{n_g}{h_g}.$
 Then $$t_g|P_c = \sigma_gt_g|U_{c,g} = \sigma_g t_g(U_{c,g}z) =
\sigma_gt_g(\frac{\text{gcd}(h_g,y)^2}{{h_g}^2}z -  \frac{u_g}{
h_g^2} \text{gcd}(h_g,y))$$
 where $\sigma_g$ is an $h_g$-th root of unity.}

\bigskip

{\em Proof.}  Since
$P_c\left(\begin{array}{ll}
     \frac{h_g}{\text{gcd}(h_g,y)}  &  \frac{u_g}{h_g}   \\
           0  &             \frac{\text{gcd}(h_g,y)}{h_g}  \\
           \end{array}
    \right) \in  n_g|h_g$ and $n_g|h_g$
 is the eigen group of $t_g(z)$
$$t_g|P_c\left(\begin{array}{ll}
     \frac{h_g}{\text{gcd}(h_g,y)}  &  \frac{u_g}{h_g}   \\
           0  &             \frac{\text{gcd}(h_g,y)}{h_g}  \\
           \end{array}
    \right) = \sigma_gt_g(z).$$
This completes the proof of the corollary.
\qed

We now consider the case that $c$ is equivalent to $\infty$
in the eigen group of $\Gamma_g$ by an Atkin-Lehner involution $W_{e_g}$.

\bigskip

{\bf  Lemma 9.}  {\em  Suppose that  $c = \frac{x}{y}$
 is equivalent to $\infty$
in the eigen group of $\Gamma_g$ by an Atkin-Lehner involution $W_{e_g}$.
 Let  an integer $u_g$ be chosen such that $e_g$ is a divisor of
  an integer $\frac{u_gy}{h_g} + \frac{zN_{\chi}\text{gcd}(h_g,y)}{h_gy_0}$
 where $$e_g = \frac{\frac{n_g}{h_g}}{
 \text{gcd}(\frac{n_g}{h_g}, \frac{y}{\text{gcd}(y,h_g)})} .$$
 Then $$P_cU_{c,g}^{-1} = P_c\left(\begin{array}{ll}
     \frac{e_gh_g}{\text{gcd}(h_g,y)}  &  \frac{u_g}{h_g}   \\
           0  &             \frac{\text{gcd}(h_g,y)}{h_g}  \\
           \end{array}
    \right) = W_{e_g} \in  n_g|h_g +e_g, f_g,  \cdot \cdot \cdot .$$
Furthermore,
 $$t_g|P_c = \sigma_g t_g(U_{c,g}z) =
 \sigma_gt_g(\frac{\text{gcd}(h_g,y)^2}{e_g{h_g}^2}z- \frac{u_g}{
e_gh_g^2} \text{gcd}(h_g,y))$$
 where $\sigma_g$ is an $h_g$-th root of unity. }

\bigskip

{\em Proof.} Let us first show that such an $u_g$ exists.  We will need
$u_g$ such that
$$ yu_g + z\frac{N_{\chi}}{y_0}\text{gcd}(h_g,y) \equiv 0 \, \pmod {e_gh_g}.$$
This follows from
$$ \frac{y}{\text{gcd}(h_g,y)}u_g + z\frac{N_{\chi}}{y_0} \equiv 0
 \pmod { \frac{h_g}{\text{gcd}(h_g,y)}e_g}.$$
Since $\text{gcd}(\frac{n_g}{h_g}, \frac{y}{\text{gcd}(h_g,y)} ) =
\frac{n_g}{h_ge_g}$
 and $e_g$ is an exact divisor of $\frac{n_g}{h_g}$, we see that
$$\text{gcd}(\frac{y}{\text{gcd}(h_g,y)},e_g) = 1.$$
Therefore $\frac{y}{\text{gcd}(h_g,y)}$ is invertible modulo
  $\frac{h_g}{\text{gcd}(h_g,y)}e_g, $ hence $u_g$ exists as required.

 The matrix
 $$P_c\left(\begin{array}{ll}
     \frac{e_gh_g}{\text{gcd}(h_g,y)}  &  \frac{u_g}{h_g}   \\
           0  &             \frac{\text{gcd}(h_g,y)}{h_g}  \\
           \end{array}
    \right)  =
  \left(\begin{array}{ll}
     \frac{xe_gh_g}{\text{gcd}(h_g,y)}  &   \frac{xu_g+w\text{gcd}(h_g,y)}{h_g}
  \\
     \frac{ye_gh_g}{\text{gcd}(h_g,y)} & \frac{yu_g}{h_g} +
                              \frac{zN_{\chi}\text{gcd}(h_g,y)}{y_0h_g} \\
           \end{array} \right)$$
 has the property

\vspace{.1cm}

(1).   $\frac{u_gy}{h_g} + \frac{zN_{\chi}\text{gcd}(h_g,y)}{h_gy_0}$
 is a multiple of $e_g$ by choice of $u_g$, and,

(2). $\frac{ye_gh_g}{\text{gcd}(h_g,y)} $ is a multiple of $n_g$, since
 $\text{gcd}(\frac{n_g}{h_g}, \frac{y}{\text{gcd}(y,h_g)}) =
 \frac{n_g}{h_ge_g}.$

\vspace{.1cm}

Therefore

$$P_c\left(\begin{array}{ll}
     \frac{e_gh_g}{\text{gcd}(h_g,y)}  &  \frac{u_g}{h_g}   \\
           0  &             \frac{\text{gcd}(h_g,y)}{h_g}  \\
           \end{array}
    \right)  = W_{e_g} \in  n_g|h_g + e_g, f_g, \cdot \cdot \cdot \, .$$
  Since $c$ is equivalent to $\infty$ in the eigen group
 of $\Gamma_g$ by $W_{e_g}$, the transformation formula follows
 easily.
\qed

\vspace{.2cm}

{\bf Remark.}  It is easy to see that Lemma 7 and Corollary 8
 are included in Lemma 9 if $e_g = 1$, in which case every element
 of $n_g|h_g$ is called an Atkin-Lehner involution for $e_g = 1.$
 This abuse of words will be used occasionally for the balance of
 the paper.

\bigskip

 The singular part of
 $ \frac{1}{|\Bbb M|} \sum_{g \in \Phi_c}
  \chi(g)(t_g|P_c)(z)$ at $z = \infty i$
  is now determined by
    $$\text{sing}_{P_c}t_{\chi}
  =  \frac{1}{|\Bbb M|} \sum_{g \in \Phi_{c}}
  \chi(g)\frac{\sigma_g}{ e^{2\pi i U_{c,g} z}}.$$

\vspace{.2cm}

We give a few examples in the calculation of
 $\text{sing}_{P_c}t_{\chi}$'s.
 Note that the first example will be used later
 in the determination
 of the invariance groups.

\bigskip

{\bf Example 1.} Suppose  $c$ = $\frac{0}{1}$. Then
$\text{sing}_{P_0}t_{\chi} =  \frac{1}{|\Bbb M|} \sum_{g \in \Phi_{c}}
  \chi(g)\frac{\sigma_g}{q^{\frac{1}{n_gh_g}}}.$

\bigskip

{\em Proof.}   In this case $y = 1$. We may choose $P_0$ =
$\left(\begin{array}{ll}
           x   &  w  \\
           y  &   \frac{xN_{\chi}}{y_0} \\
    \end{array}
     \right)$  =
$\left(\begin{array}{rr}
           0   &  -1  \\
           1  &   0 \\
    \end{array}
  \right)$
 and  $u_g = 0.$
 If the condition of Corollary 8 holds, then
 $n_g = h_g$ and
$$t_g|P_0 = \sigma_g t_g(\frac{z}{h_g^2}) =
  \sigma_g t_g(\frac{z}{n_gh_g}). $$
  The only $ g \in \Phi_0 \subseteq \Bbb M$ satisfying
 $n_g = h_g$ are $1A$ and  $3C$. We have
$$t_{1A}|P_0 = t_{1A} \quad \text{and} \quad t_{3C}|P_0 =
 \sigma_{3C} t_{3C}(\frac{1}{9}z).$$
On the other hand, if the condition of Lemma 9 holds, then
 $ n_g = e_gh_g$ and
$$ t_g|P_0 = \sigma_g t_g(\frac{z}{e_gh_g^2})
 = \sigma_g t_g(\frac{z}{n_gh_g}).$$
Note that for all the
 remaining $g \in \Phi_0 \setminus \{ 1A, 3C \}
\subseteq \Bbb M $, $0$ is equivalent
 to $\infty$ in $\Gamma_g$ by the  Atkin-Lehner
 involution $W_{\frac{n_g}{h_g}}$ and hence
 the condition of Lemma 9 holds.
\qed

\bigskip

{\bf Remark.}  Example 1 shows that $0$ is a pole of the
 McKay-Thompson series $t_{\chi}(z)$ for every $\chi$, since $n_gh_g \ne 1$
 for $g \ne 1 $ and so the coefficient of $q^{-1}$ is nonzero.

\bigskip

{\bf Example 2.}  Let $\chi$ be the trivial character and
 let $c = \frac{1}{3}$.  Then $P_{\frac{1}{3}}$ =
 $\left(\begin{array}{ll}
           1   &  [\frac{N_0}{81}] \\
           3  &   \frac{N_0}{27} \\
    \end{array}
  \right)$ ($[x]$ is the integral part of $x$ and $\frac{N_0}{27}
 \equiv 1$ ( mod $3$) ) and
 $t_{84A}|P_{\frac{1}{3}}$ =
 $\sigma_{84A}  t_{84A}(\frac{1}{56}z + \frac{1}{2}).$

\bigskip

{\em Proof.}  We know $\Gamma_{84A} < 84|2+$.
Since $\text{gcd}(\frac{84}{2}, \frac{3}{\text{gcd}(3,2)})$ = $3$
 = $\frac{84}{2e}$,
  we  see  that $e$ is $14$, and can choose  $u_g = -28.$
 $$P_{\frac{1}{3}}\left(\begin{array}{rr}
           28   &  -\frac{28}{2}  \\
           0  &   \frac{1}{2} \\
    \end{array}
  \right) = W_{14} \in  84|2+.$$
 The rest follows easily.
\qed

\bigskip

{\bf Remark.}
(1).   The invariance group $\Gamma_g$ of the harmonics
 $n|h + e,f, \cdot \cdot \cdot $ are
 not fully determined.
(For each $g$, one can write down  a set of generators
 of  the invariance group  $\Gamma_g$ easily.
  But determining whether or not a given element in $SL_2(\Bbb R)$
 is a word of those generators
 is nontrivial.)
  Hence we have to settle for $\sigma_g$
 being an $h_g$-th root of unity.

(2).  $\sigma_g = 1$ if  $h_g$ = $1$.

(3).  Let $p \in \{11,17,19,23,29,31,41,47,59,71\}. $
 Applying Lemma 5, one can prove
$\Phi_0 =\Phi_{\frac{x}{p}}$
  and $\Phi_{\frac{1}{32}}=\Phi_{\frac{1}{64}}$.

\bigskip

\section{invariance group}

\bigskip

Let $f$ be a modular function and let
 $K_f$ be a subgroup of the invariance group (in $SL_2(\Bbb R)$)
  $\Gamma_f$ of $f$ of finite
 index.
 We shall determine   the invariance group  as follows.
   Define

\vspace{.2cm}

$C_f$ = the set of all  cusps of $K_f$, and,

$C_{0}$ = $\{ c \in C_f| c $ is a
 pole of  $  f \}.$

\bigskip

{\bf Lemma 10.}  {\em We have $\Gamma_f C_f \subseteq C_f$  and
   $ \Gamma_f C_{0} \subseteq
 C_{0}.$}

\bigskip

{\em Proof.}  Since $[\Gamma_f: K_f] < \infty,$ $C_f$ is also the
 set of all cusps of $\Gamma_f.$  The second statement is obvious.
\qed

\bigskip

{\bf Lemma 11.} {\em Let $f$ be a modular function and let $\Gamma_f$
 be its invariance group.  Suppose that $K_f \le
 \Gamma_f$. Let $\alpha = \frac{a_1}{c_1}$} ($a_1 \ne 0$)
  {\em and $\beta = \frac{a_2}{c_2}$
 be two inequivalent cusps of $K_f$.  Let
$$M_1= \left(\begin{array}{ll}
                 a_1  & b_1 \\
                 c_1 &  d_1 \\
             \end{array}
             \right) \quad and \quad
 M_2= \left(\begin{array}{ll}
                 a_2  & b_2 \\
                 c_2 &  d_2 \\
             \end{array}
             \right).$$
Then $\frac{a_1}{c_1}$ and $\frac{a_2}{c_2}$ are equivalent
 with respect to $\Gamma_f$  if and only if
 the $q$-expansion of
$f|_{M_1}$ is derived from  that of $f|_{M_2}$ under
 the substitution $ z \to az+b$ for some  numbers
 $a$ and $b$} ({\em if $c_i = 0,$ then $a_i = 1$ and $\frac{a_i}{c_i}
 = \infty$}).

\bigskip

{\em Proof.}  Let $A \in \Gamma_f$ be such that $A\alpha = \beta$.
 Define  the  matrix $B$ such that
 $$A = M_2 \left(\begin{array}{rr}
                     1 & 0 \\
                 -\frac{c_1}{a_1} &  1 \\
                 \end{array}
                 \right)B.$$
Since   $A\alpha = \beta,$  it follows that $B\alpha = \alpha.$  Hence
$$B  = \left(\begin{array}{ll}
                     1 & 0 \\
                  \frac{c_1}{a_1} &  1 \\
                 \end{array}
                 \right)
          \left(\begin{array}{ll}
                    m_{11} & m_{12} \\
                    0   &  m_{22} \\
                 \end{array}
                 \right)
          \left(\begin{array}{rr}
                     1 & 0 \\
                  -\frac{c_1}{a_1} &  1 \\
                 \end{array}
                 \right)$$
 for some  $m_{11}$, $m_{12}$ and $m_{22}$.
 In particular,
 $$A = M_2\left(\begin{array}{ll}
                     m_{11} & m_{12} \\
                        0 &  m_{22} \\
                 \end{array}
                 \right)
            \left(\begin{array}{rr}
                     1 & 0 \\
                  -\frac{c_1}{a_1} &  1 \\
                 \end{array}
                 \right)$$
and
  $$A\alpha = M_2\left(\begin{array}{ll}
                     m_{11} & m_{12} \\
                     0 & m_{22} \\
                 \end{array}
                 \right)\infty.$$
It follows that
$f|_{M_1}$ = $f|_{AM_1}$ =  $$f|{M_2
              \left(\begin{array}{rr}
                     1 & 0 \\
                  -\frac{c_1}{a_1} &  1 \\
                 \end{array}
                 \right)BM_1}
           =   f|{M_2\left(\begin{array}{ll}
                     m_{11} & m_{12} \\
                       0 &  m_{22} \\
                 \end{array}
                 \right)
                \left(\begin{array}{ll}
                     a_1 & b_1 \\
                       0 & a' \\
                 \end{array}
                 \right)}
           = f|{M_2\left(\begin{array}{ll}
                     m_{11}' & m_{12}' \\
                      0 & m_{22}'     \\
                    \end{array}
                    \right)},$$
               for some  $b_1$, $a'$, $m_{11}'$, $m_{12}'$
 and $m_{22}'$. Consequently,
 $f|_{M_1}$ is derived from  that of $f|_{M_2}$ under
 the substitution $ z \to az+b$
 where $a = \frac{m_{11}'}{m_{22}'}$ and $b =
  \frac{m_{12}'}{m_{22}'}.$  Conversely, one sees easily that
$\alpha$ and $\beta$ is equivalent to each other by
$$M_1\left(\begin{array}{rr}
                     \frac{1}{a} & -\frac{b}{a} \\
                      0 & 1     \\
                    \end{array}
                    \right)M_2^{-1} \in \Gamma_f.$$
\qed

\bigskip

The invariance group $\Gamma_f$   can now be determined as
 follows :

\vspace{.2cm}

(1).  Determine $C_f$, the set of all cusps of $K_f$.

(2).  Determine the subset $C_{0}.$

(3).  Determine the $q$-expansion of $f$ at all $c_i$  in $C_{0}$
 by suitable matrices $M_i$ such that $M_i\infty = c_i.$

(4).  Apply Lemma 11 and  determine the set $E_{0} =
 \{ c \in C_0 |$ the $q$-expansion of $f$ at $c$ is derived from
  the $q$-expansion of $f$ at $0$ under the
 substitution $z \to az + b \}$ and  the set $ A_{0}$ =
 $\{A_c \in \Gamma_f, A_c0 = c\}.$  Note that it is sufficient
 to determine at most one matrix $A_c$  for each
 representative of inequivalent cusps.

(5).  Determine $(\Gamma_f)_{0}$ =
 $\langle B| B = \left(\begin{array}{ll}
                          1 & 0 \\
                          m & 1 \\
                 \end{array} \right), m \in r\Bbb Z,
 \text{ for some (fixed) }  r \in \Bbb Q\rangle.$
Note that this can be
 achieved
 by investigating the $q$-expansion of $f$ at
 $0$.  Note also that $B$ is of the form given since $\Gamma_f$
 is discrete.

\bigskip

{\bf Remark.}
(1).  The McKay-Thompsom series
 $t_{\chi}(z)$ has a pole at $0$ for every $\chi$
 as stated in the remark
 right after Example 1.

(2). Since Lemma 11
 applies only when one of the cusp is nonzero,
   one can not take $c$ to be $0$ in (4) above.

(3). One can replace $0$ by any cusp and apply our procedure
 to find the invariance groups.

\bigskip

{\bf Lemma 12. } {\em  $\Gamma_f = \langle K_f, B, A_c, c \in E_{0}
 \rangle.$ }

\bigskip

{\em Proof.}  For  any $\sigma \in \Gamma_f \setminus
 \langle K_f, A_c, c \in E_{0}\rangle$.
 Applying Lemma 10, $\sigma0$ is again a cusp.  Hence $\sigma0$
 must  be $\langle K_f, A_c, c \in E_{0}\rangle$-equivalent to $0.$
 Choose $\delta \in \langle K_f, A_c, c \in E_{0}\rangle$
 such that $\delta\sigma 0 = 0.$
  Then $\delta\sigma \in (\Gamma_f)_0$.
 Hence $\Gamma_f = \langle K_f, B, A_c, c \in E_{0}
 \rangle$ holds.
\qed

\bigskip

{\bf Theorem 13.} (Helling's Theorem [5])
  {\em  The maximal discrete groups of $PSL(2,\Bbb C)$
 commensurable with the modular group $SL(2,\Bbb Z)$ are just the
  images of the conjugates of $\Gamma_0(N)+$ for square free $N.$ }

\bigskip

{\bf Corollary 14.} {\em  For each irreducible character $\chi$
 of $\Bbb M$,
 the set of  prime divisors of the index
$[\Gamma_{\chi}:\Gamma_0(N_{\chi})]$
 is a subset of $\{2,3,5,7 \}$. }

\bigskip

{\em Proof.}  By Helling's Theorem, any  maximal subgroup that
contains $\Gamma_{\chi}$ is a conjugate of some $\Gamma_0(n)+$.
        Conway has shown in [1] that  $n$ must be a divisor of $N_{\chi}$.
Now compare the volumes of
 the fundamental domains  of  $SL_2(\Bbb Z)$, $\Gamma_0(n)$,
 and $\Gamma_0(N_{\chi})$.  Noting that the conjugation does not
 change the volume and that the normalizer of $\Gamma_0(n)$
 changes the volume of the fundamental domain by a
 factor involving only primes $2$ and $3$, we obtain our
 lemma since the index $[SL_2(\Bbb Z): \Gamma_0(N_{\chi})]$
 involves only primes $2$, $3$, $5$, and $7$.
\qed

\bigskip

 Let $\chi$ be an irreducible character
  of $\Bbb M$ and  $\Gamma_{\chi}$
 be the invariance group of $t_{\chi}(z)$.
 We are now ready to prove :

\vspace{.2cm}

(1).  $A_0 = \emptyset,$ and,

(2).  $(\Gamma_{\chi})_{0} = (\Gamma_0(N_{\chi}))_{0}.$

\bigskip

{\bf Lemma 15.}  {\em Let $\chi$ be an irreducible
 character of $\Bbb M$
 and let $c$ be a cusp of $\Gamma_0(N_{\chi}),$ not equivalent to $0$.
 Then $A_0 = \{ A_c | A_c0 = c\} = \emptyset.$ }

\bigskip

{\em Proof.}  We first recall that the singular part of
 $ \frac{1}{|\Bbb M|} \sum_{g \in \Phi_c}
  \chi(g)(t_g|P_c)(z)$ at $z = \infty i$
  is given  by
    $$\text{sing}_{P_c}t_{\chi}
  =  \frac{1}{|\Bbb M|} \sum_{g \in \Phi_{c}}
  \chi(g)\frac{\sigma_g}{ e^{2\pi i U_{c,g} z}}.$$
 Applying Lemmas 7, 9 and 10 and Corollary 8,
 we see that  it suffices to
 show that $\text{sing}_{P_0}t_{\chi}$ can not be derived from
 $\text{sing}_{P_c}t_{\chi}$
 under the substitution $ z \rightarrow az + b$ if $c \ne 0$.
  This is achieved
 by a $case$-$by$-$case$ study.
  We give an example to indicate how  the lemma is proved.

\vspace{.2cm}

{\bf Example 3.}  Let $\chi$ be the trivial
 character of $\Bbb M$.  Then $A_0 = \emptyset.$

\vspace{.2cm}

{\em Proof.}  We first note that for any irreducible character
 $\chi$ of $\Bbb M$ and $c \in \Bbb Q \cup \{ \infty \},$  we have :

(1). The lowest terms in $\text{sing}_{P_c}t_{\chi}$ and
   $\text{sing}_{P_0}t_{\chi}$
 are  of the form $\frac{r}{q}$ for some number $r
 \in \Bbb Q, $ and

(2). Terms in
   $\text{sing}_{P_c}t_{\chi}$ and
 $\text{sing}_{P_0}t_{\chi}$ are all of the form
 $\frac{r}{q^{\frac{1}{t}}}$ (Corollary 8, Lemma 9),
 for some $r \in \Bbb Q$, and $t \in \Bbb N$.

 Since  the lowest term in  $\text{sing}_{P_0}t_{\chi}$ is
 $ \frac{r}{q},$ $r \ne 0,$
 the transformation that sends   $\text{sing}_{P_c}t_{\chi}$ to
  $\text{sing}_{P_0}t_{\chi}$
 is of the form $ z \rightarrow az+ b$  where $a$ is some positive
  integer.

Let $\chi$ is the trivial character, then by Corollary 8 and Lemma 9,
 one has

 $$\text{sing}_{P_0}t_{\chi}
  =  \frac{1}{|\Bbb M|} \sum_{g \in \Phi_{0}}
  \frac{\sigma_g}{q^{\frac{1}{n_gh_g}   }} =
  \frac{1}{|\Bbb M|}(\frac{1}{q} + \frac{2|\Bbb M|}{|C_{\Bbb M}(71A)|}
  \frac{\sigma_{71A}}{q^{\frac{1}{71}}} + \cdot \cdot \cdot ), $$

 and for  any cusp $c = \frac{x}{y}$,  $\text{gcd}(x,y) =1,$
  $$\text{sing}_{P_c}t_{\chi}
  =  \frac{1}{|\Bbb M|} \sum_{g \in \Phi_{c}}
  \frac{\sigma_g}{ e^{2\pi i U_{c,g}z}}.$$

Suppose  $\text{sing}_{P_0}t_{\chi}$ can be derived from
  $\text{sing}_{P_c}t_{\chi}$  under   $ z \rightarrow az + b$.

(3).  We will show  gcd($y$, $71$) = 1.  Suppose false. Then
 $71|y$ and $71A,$ $71B \in \Phi_c.$
Since  $\frac{1}{q^{\frac{1}{71}}}$ appears in
 $\text{sing}_{P_0}t_{\chi}$, there exists, by Lemma 11,
 some $g$ in  $\Phi_c$ such that
  $$\frac{\sigma_g}{e^{2\pi i U_{c,g}z}}|(z \rightarrow az + b)
 = \frac{r}{q^{\frac{1}{71}}}  \eqno(*) $$
where $r$ is some constant. Hence
$$\frac{ \text{gcd}(h_g,y)^2a}{e_gh_g^2} = \frac{1}{71}$$
where
 $$e_g = \frac{\frac{n_g}{h_g}}{
 \text{gcd}(\frac{n_g}{h_g}, \frac{y}{\text{gcd}(y,h_g)})} .$$
Since $a$ is an integer, this implies that
$ g = 71A$ or $ 71B. $
 Since $\Gamma_{71A} = \Gamma_{71B} = 71+$, we have $ n_{71A}$
 = $n_{71B}$ = $71$, $h_{71A} $ = $h_{71B}$ = $1$,
  $e_{71A}$ = $e_{71B}$ = $71$, and $a = 1.$
 On the other hand, Corollary 8 implies $t_g|P_c = \sigma_g
t_g(z - u_g)$ for $g = 71A$ or $71B$.  Hence the
 transformation $(*)$ can not be done.     This forces
$$ \text{gcd}(y,71) = 1.$$

(4).  Since $\frac{1}{q^{\frac{1}{t}}}$,
  $ t \in \{ 29, 41, 59, 92, 93, 94, 95, 104, 110, 119 \}$
  all appear with nonzero
 coefficients in  $\text{sing}_{P_0}t_{\chi}$,
  we may similarly conclude that $\text{gcd}(y, p) = 1$
 for the  other prime divisors  of $N_0$.

 Hence  $\text{gcd}(y,N_0) = 1$ and  $c$ is equivalent to $0$.
    Consequently, $A_0 = \emptyset.$

\bigskip

{\bf Lemma 16.}  {\em $(\Gamma_{\chi})_{0}
 = (\Gamma_0(N_{\chi})_{0}.$ }

\bigskip

{\em Proof.}  Suppose not.  Applying Corollary 14, we see that
$ (\Gamma_{\chi})_0$ contains $B_r =
 \left( \begin{array}{ll}
               1   &   0 \\
           \frac{N_{\chi}}{r}  &  1 \\
          \end{array}
    \right)$ for $r$ =  $2$, $3$, $5$ or $7$.
  This implies that the cusp $\infty$ is equivalent to
 $\frac{r}{N_{\chi}}$ in $\Gamma_{\chi}$.  Therefore
$\text{sing}_{\infty}t_{\chi}$
 must be derived from $\text{sing}_{P_{\frac{r}{N_{\chi}}}}t_{\chi}$
 under the substitution $ z \rightarrow az+b$.
  We can now apply an analogous procedure (using $y = \frac{r}{N_{\chi}}$)
 as in Example 3 to get a contradiction.
\qed

\bigskip

{\bf Remark.}  One can also prove Lemma 16 by claiming
 that $B_r$ does not leave $t_{\chi}(z)$ invariant. Note that
 it is easy to show the claim  since $B_r$ leaves
 most of the $t_g(z)$ invariant except for those $g$'s
 such that $n_gh_g$ is not a divisor of $\frac{N_{\chi}}{r}.$

\bigskip

Combining Lemma 15 and 16, we have :

\bigskip

{\bf Theorem 17.} {\em Let $\chi$ be an irreducible character of $\Bbb M$.
   Then $\Gamma_{\chi} = \Gamma_0(N_{\chi}).$ }

\bigskip

$N_{\chi}$ can be found in $Table$ $1$.

\bigskip

{\bf Remark.}  (1). In Lemma 15 and 16, $0$ is a better choice than the other
 cusps ($\infty$, for example)
 since among all the $\text{sing}_{P_c}t_{\chi}$'s,
 $\text{sing}_{P_0}t_{\chi}$
 is the one that involves most  nonzero terms.

(2).  $N_{\chi}$ and its prime decomposition is calculated by a software
 called GAP.

\bigskip
\bigskip
\bigskip
\bigskip

\newpage

$$Table \quad 1$$

$\begin{array}{lll}

\chi_i  &  N_{\chi_i}  &  N_{\chi_i} \, (prime \, decomposition) \\
        &          &                                \\
1 &  2331309585756753201600 &   2^63^3 5^2 7. 11. 13. 17. 19. 23. 29. 31. 41.
47.
  59. 71 \\
2 & 11841091337275200    & 2^63^3 5^2 7. 11. 13. 17. 19. 23. 29. 31. 41 \\
3 & 437868837806400 &      2^63^3 5^2 7. 11. 13. 17. 19. 23. 29. 47 \\
4 & 467584848090400 & 2^5 3^3 5^2 7. 11. 13. 17. 19. 23. 41. 71 \\
5 & 38732026132800 & 2^63^3 5^2 7. 11. 13. 17. 19. 47. 59 \\
6 &  20350725595200 & 2^63^3 5^2 7. 11. 13. 17. 19. 31. 47 \\
7 & 87358471200 &  2^5 3^2 5^2 7. 11. 13. 17. 23. 31 \\
8 &  7820482269600 &2^5 3^2 5^2 7. 11. 13. 17. 29. 31. 71 \\
9 & 18526958049600&  2^63^2 5^2 7. 11. 13. 23. 29. 41. 47 \\
10 &222987885120& 2^63^3 5. 7. 11. 13. 19. 23. 59 \\
11 &8490081600& 2^63^2 5^2 7. 11. 13. 19. 31 \\
12 &19445025600& 2^63^2 5^2 7. 11. 13. 19. 71 \\
13 &9958865716800& 2^63^3 5^2 7. 11. 13. 17. 19. 23. 31 \\
14 &73513400& 2^5 3^3 5. 7. 11. 13. 17 \\
15 &2244077793757800& 2^3 3^3 5^2 7. 11. 13. 17. 19. 23. 29. 41. 47 \\
16 &3749442460305984& 2^63^3 13. 23. 29. 31. 41. 47. 59. 71 \\
17 &3749442460305984& 2^63^3 13. 23. 29. 31. 41. 47. 59. 71 \\
18 &726818400& 2^5 3^3 5^2 7. 11. 19. 23 \\
19 &9182927033280& 2^63^3 5. 7. 11. 13. 19. 29. 41. 47 \\
20 &35703027360& 2^5 3^2 5. 7. 11. 13. 17. 31. 47 \\
21 &98066928960& 2^63^3 5. 7. 11. 13. 17. 23. 29 \\
22 & 22789166400& 2^63^3 5^2 7. 11. 13. 17. 31 \\
23 &451392480& 2^5 3^3 5. 7. 11. 23. 59 \\
24 &295495200& 2^5 3^2 5^2 7. 11. 13. 41 \\
25 &253955520& 2^63^3 5. 7. 13. 17. 19 \\
26 &479256378753600& 2^63^2 5^2 7. 19. 31. 41. 47. 59. 71 \\
27 &479256378753600& 2^63^2 5^2 7. 19. 31. 41. 47. 59. 71 \\
28 &27003936960& 2^63^2 5. 7. 11. 13. 17. 19. 29 \\
29 &81995760& 2^4 3^3 5. 7. 11. 17. 29 \\
30 &69618669120& 2^63^2 5. 7. 11. 13. 19. 31. 41 \\
31 &21416915520& 2^63^2 5. 7. 11. 13. 17. 19. 23 \\
32 &214885440& 2^63^3 5. 7. 11. 17. 19 \\
33 &2882880& 2^63^2 5. 7. 11. 13 \\
34 &332640& 2^5 3^3 5. 7. 11 \\
35 &786240& 2^63^3 5. 7. 13 \\
36 &11147099040& 2^5 3^3 5. 7. 11. 23. 31. 47 \\
37 &331962190560& 2^5 3^2 5. 7. 11. 13. 17. 19. 23. 31 \\
38 &333637920& 2^5 3^3 5. 7. 11. 17. 59 \\
39 &845013600& 2^5 3^3 5^2 19. 29. 71 \\
40 &845013600& 2^5 3^3 5^2 19. 29. 71 \\
41 &16676856385200& 2^4 3^3 5^2 11. 23. 31. 47. 59. 71 \\
42 &16676856385200& 2^4 3^3 5^2 11. 23. 31. 47. 59. 71 \\
43 &186902100& 2^2 3^3 5^2 7. 11. 29. 31 \\
44 &46955594400& 2^5 3. 5^2 11. 13. 41. 47. 71 \\
\end{array}$

\newpage

$\begin{array}{lll}
45 &46955594400& 2^5 3. 5^2 11. 13. 41. 47. 71 \\
46 &54880846020& 2^2 3^2 5. 7. 11. 13. 17. 19. 23. 41 \\
47 &105386400& 2^5 3^3 5^2 7. 17. 41 \\
48 &105386400& 2^5 3^3 5^2 7. 17. 41 \\
49 &49584815280& 2^4 3^3 5. 7. 11. 13. 17. 19. 71 \\
50 &6404580& 2^2 3^2 5. 7. 13. 17. 23 \\
51 &12916800& 2^63^3 5^2 13. 23 \\
52 &646027200& 2^63^2 5^2 7. 13. 17. 29 \\
53 &228731328& 2^63^2 7. 17. 47. 71 \\
54 &228731328& 2^63^2 7. 17. 47. 71 \\
55 &19044013248& 2^63^3 13. 23. 29. 31. 41 \\
56 &10944013248& 2^63^3 13. 23. 29. 31. 41 \\
57 &25077360& 2^4 3. 5. 7. 11. 23. 59 \\
58 &198918720& 2^63^3 5. 7. 11. 13. 23 \\
59 &19433872080& 2^4 3^3 5. 7. 23. 29. 41. 47 \\
60 &19433872080& 2^4 3^3 5. 7. 23. 29. 41. 47 \\
61 &2784600& 2^3 3^2 5^2 7. 13. 17 \\
62 &245044800& 2^63^2 5^2 7. 11. 13. 17 \\
63 &57266969760& 2^5 3^3 5. 7. 11. 13. 17. 19. 41 \\
64 &157477320& 2^3 3^2 5. 7. 11. 13. 19. 23 \\
65 &818809200& 2^4 3^2 5^2 11. 23. 29. 31 \\
66 &263877213600& 2^5 3^2 5^2 7. 11. 13. 19. 41. 47 \\
67 &1588466880& 2^63^2 5. 7. 11. 13. 19. 29 \\
68 &33005280& 2^5 3. 5. 7. 11. 19. 47 \\
69 &937440& 2^5 3^3 5. 7. 31 \\
70 &32864832& 2^63^3 7. 11. 13. 19 \\
71 &182584514112& 2^63. 7. 11. 13. 17. 29. 41. 47 \\
72 &182584514112& 2^63. 7. 11. 13. 17. 29. 41. 47 \\
73 &982080& 2^63^2 5. 11. 31 \\
74 &33542208& 2^63^3 7. 47. 59 \\
75 &33542208& 2^63^3 7. 47. 59 \\
76 &7650720& 2^5 3^3 5. 7. 11. 23 \\
77 &931170240& 2^63^2 5. 7. 11. 13. 17. 19 \\
78 &33754921200& 2^4 3^2 5^2 7. 11. 13. 17. 19. 29 \\
79 &42325920& 2^5 3^2 5. 7. 13. 17. 19 \\
80 &4969440& 2^5 3^2 5. 7. 17. 29 \\
81 &63126554400& 2^5 3^2 5^2 7. 13. 23. 59. 71 \\
82 &63126554400& 2^5 3^2 5^2 7. 13. 23. 59. 71 \\
83 &208304928& 2^5 3^2 13. 23. 41. 59 \\
84 &208304928& 2^5 3^2 13. 23. 41. 59 \\
85 &704223936& 2^63^3 13. 23. 29. 47 \\
86 &704223936& 2^63^3 13. 23. 29. 47 \\
87 &1235520& 2^63^3 5. 11. 13 \\
88 &3967200& 2^5 3^2 5^2 19. 29 \\
89 &11737440& 2^5 3^3 5. 11. 13. 19 \\
90 &11737440& 2^5 3^3 5. 11. 13. 19 \\
91 &2542811040& 2^5 3^2 5. 7. 11. 17. 19. 71 \\
92 &22102080& 2^63. 5. 7. 11. 13. 23 \\
93 &1441440& 2^5 3^2 5. 7. 11. 13 \\
94 &879840& 2^5 3^2 5. 13. 47 \\

\end{array}$

\newpage

$\begin{array}{lll}

95 &16576560& 2^4 3^2 5. 7. 11. 13. 23 \\
96 &21677040& 2^4 3^2 5. 7. 11. 17. 23 \\
97 &5267201940& 2^2 3^2 5. 7. 11. 13. 23. 31. 41 \\
98 &7900200& 2^3 3^3 5^2 7. 11. 19 \\
99 &1660401600& 2^63. 5^2 11. 13. 41. 59 \\
100 &1660401600& 2^63. 5^2 11. 13. 41. 59 \\
101 &932769600& 2^63. 5^2 7. 17. 23. 71 \\
102 &7601451872175& 3^3 5^2 7. 17. 19. 29. 41. 59. 71 \\
103 &7601451872175& 3^3 5^2 7. 17. 19. 29. 41. 59. 71 \\
104 &6511680& 2^63^2 5. 7. 17. 19 \\
105 &5844589984800& 2^5 3^2 5^2 7. 19. 31. 47. 59. 71 \\
106 &5844589984800& 2^5 3^2 5^2 7. 19. 31. 47. 59. 71 \\
107 &2434219200& 2^63^2 5^2 7. 19. 31. 41 \\
108 &2434219200& 2^63^2 5^2 7. 19. 31. 41 \\
109 &280800& 2^5 3^3 5^2 13 \\
110 &947520& 2^63^2 5. 7. 47 \\
111 &1016747424& 2^5 3^3 7. 11. 17. 29. 31 \\
112 &386100& 2^2 3^3 5^2 11. 13 \\
113 &6568800& 2^5 3. 5^2 7. 17. 23 \\
114 &1374912& 2^63^2 7. 11. 31 \\
115 &151200& 2^5 3^3 5^2 7 \\
116 &92400& 2^4 3. 5^2 7. 11 \\
117 &411840& 2^63^2 5. 11. 13 \\
118 &19562400& 2^5 3^2 5^2 11. 13. 19 \\
119 &12524852340& 2^2 3^3 5. 7. 11. 13. 17. 29. 47 \\
120 &41801760& 2^5 3^2 5. 7. 11. 13. 29 \\
121 &75698280& 2^3 3^3 5. 7. 17. 19. 31 \\
122 &8148853440& 2^63^2 5. 7. 13. 17. 31. 59 \\
123 &864175548600& 2^3 3^3 5^2 7. 13. 17. 31. 47. 71 \\
124 &3456702194400& 2^5 3^3 5^2 7. 13. 17. 31. 47. 71 \\
125 &3456702194400& 2^5 3^3 5^2 7. 13. 17. 31. 47. 71 \\
126 &119700& 2^2 3^2 5^2 7. 19 \\
127 &13695552& 2^63^2 13. 31. 59 \\
128 &752016096& 2^5 3^3 13. 23. 41. 71 \\
129 &752016096& 2^5 3^3 13. 23. 41. 71 \\
130 &19320840& 2^3 3^2 5. 7. 11. 17. 41 \\
131 &1004683680& 2^5 3^2 5. 7. 11. 13. 17. 41 \\
132 &164160& 2^63^3 5. 19 \\
133 &14379596431200& 2^5 3^3 5^2 7. 11. 13. 17. 19. 29. 71 \\
134 &14208480& 2^5 3^3 5. 11. 13. 23 \\
135 &497653200& 2^4 3^3 5^2 11. 59. 71 \\
136 &497653200& 2^4 3^3 5^2 11. 59. 71 \\
137 &995276700& 2^2 3^3 5^2 11. 23. 31. 47 \\
138 &5109369408& 2^63^3 11. 13. 23. 29. 31 \\
139 &514080& 2^5 3^3 5. 7. 17 \\
140 &59017104226080& 2^5 3^3 5. 7. 11. 13. 17. 19. 29. 31. 47 \\
141 &1151710560& 2^5 3^2 5. 7. 11. 13. 17. 47 \\
142 &3767400& 2^3 3^2 5^2 7. 13. 23 \\
143 &7600320& 2^63^2 5. 7. 13. 29 \\
144 &11823840& 2^5 3^3 5. 7. 17. 23 \\

\end{array}$

\newpage

$\begin{array}{lll}
145 &8558550& 2. 3^2 5^2 7. 11. 13. 19 \\
146 &664020& 2^2 3^2 5. 7. 17. 31 \\
147 &4320& 2^5 3^3 5 \\
148 &655188534& 2. 3^3 7. 11. 13. 17. 23. 31 \\
149 &102240& 2^5 3^2 5. 71 \\
150 &157248& 2^63^3 7. 13 \\
151 &26489342880& 2^5 3^2 5. 7. 11. 13. 17. 23. 47 \\
152 &93276960& 2^5 3. 5. 7. 17. 23. 71 \\
153 &13137600& 2^63. 5^2 7. 17. 23 \\
154 &428400& 2^4 3^2 5^2 7. 17 \\
155 &18221280& 2^5 3. 5. 7. 11. 17. 29 \\
156 &190072512& 2^63^2 7. 17. 47. 59 \\
157 &21176100& 2^2 3^3 5^2 11. 23. 31 \\
158 &37130940& 2^2 3^3 5. 7. 11. 19. 47 \\
159 &852390& 2. 3^3 5. 7. 11. 41 \\
160 &184363200& 2^63^2 5^2 7. 31. 59 \\
161 &108803771818560& 2^63^2 5. 7. 13. 17. 19. 23. 29. 41. 47 \\
162 &1657656& 2^3 3^2 7. 11. 13. 23 \\
163 &345290400& 2^5 3^2 5^2 7. 13. 17. 31 \\
164 &90014400& 2^63^2 5^2 7. 19. 47 \\
165 &30240& 2^5 3^3 5. 7 \\
166 &4032& 2^63^2 7 \\
167 &4062240& 2^5 3^2 5. 7. 13. 31 \\
168 &3204801600& 2^63. 5^2 7. 11. 13. 23. 29 \\
169 &24196995900& 2^2 3^2 5^2 7. 11. 17. 19. 23. 47 \\
170 &668304& 2^4 3^3 7. 13. 17 \\
171 &73920& 2^63. 5. 7. 11 \\
172 &6983776800& 2^5 3^3 5^2 7. 11. 13. 17. 19 \\
173 &32959080& 2^3 3^2 5. 7. 11. 29. 41 \\
174 &3115200& 2^63. 5^2 11. 59 \\
175 &48163383908640& 2^5 3^2 5. 7. 11. 13. 17. 19. 31. 47. 71 \\
176 &427211200& 2^65^2 13. 19. 23. 47 \\
177 &858816& 2^63^3 7. 71 \\
178 &21416915520& 2^63^2 5. 7. 11. 13. 17. 19. 23 \\
179 &14400& 2^63^2 5^2 \\
180 &14400& 2^63^2 5^2 \\
181 &154881891350& 2. 3^3 5^2 11. 13. 19. 29. 31. 47 \\
182 &2009280& 2^63. 5. 7. 13. 23 \\
183 &6339168& 2^5 3^3 11. 23. 29 \\
184 &26429760& 2^63^3 5. 7. 19. 23 \\
185 &32730048& 2^63^3 13. 31. 47 \\
186 &7425600& 2^63. 5^2 7. 13. 17 \\
187 &8237275200& 2^63^2 5^2 7. 11. 17. 19. 23 \\
188 &15120& 2^4 3^3 5. 7 \\
189 &54774720& 2^63^2 5. 7. 11. 13. 19 \\
190 &27989280& 2^5 3^3 5. 11. 19. 31 \\
191 &34272& 2^5 3^2 7. 17 \\
192 &3500640& 2^5 3^2 5. 11. 13. 17 \\
193 &1049200425& 3^3 5^2 7. 13. 19. 29. 31 \\
194 &1404480& 2^63. 5. 7. 11. 19 \\
\end{array}$

\bigskip

{\small
\noindent DEPARTMENT OF MATHEMATICS,
 THE OHIO STATE UNIVERSITY,
 COLUMBUS, OHIO  43210,  USA\\}
\noindent {\em E-mail address} : haradako@@math.ohio-state.edu \\

{\small
\noindent DEPARTMENT OF MATHEMATICS,
NITIONAL UNIVERSITY OF SINGAPORE,
SINGAPORE 0511,
REPUBLIC OF SINGAPORE\\}
\noindent {\em E-mail address} : matlml@@leonis.nus.sg


\newpage


\begin{thebibliography}{99}

\bigskip

\bigskip
\baselineskip=12pt

\bibitem{C}   J.~H. Conway, {\em Understanding Groups like $\Gamma_0(N)$},
  (preprint).

\bibitem{CN}  J.~H. Conway and S.~P. Norton, {\em Monstrous Moonshine,}
  Bull. London Math. Soc. {\bf 11} (1979), 308-338.

\bibitem{FLM} I.~Frenkel, J.~Lepowsky and A.~ Meurman,
 {\em Vertex Operator Algebras and the Monster},  Academic Press, Inc.
  (1988).

\bibitem{HA} K.~Harada, {\em Modular Functions, Modular Forms and
Finite Group},  Lecture Notes at The Ohio State University, (1987).

\bibitem{H}  H.~Helling,  {\em On the Commensurablility Classes
 of Rational Modular Group}, J. London Math. Soc.
 {\bf 2} (1970), 67-72.

\bibitem{MS}  J.~McKay and H.~Strauss, {\em The q-series of
 Monstrous Moonshine and the Decomposition of the Head Characters},
  Comm. Alg. {\bf 18(1)} (1990), 253-278.



\end{thebibliography}
\end{document}